\documentclass[sn-mathphys,Numbered]{sn-jnl}

\usepackage{graphicx}%
\usepackage{multirow}%
\usepackage{amsmath,amssymb,amsfonts}%
\usepackage{amsthm}%
\usepackage{mathrsfs}%
\usepackage[title]{appendix}%
\usepackage{xcolor}%
\usepackage{textcomp}%
\usepackage{manyfoot}%
\usepackage{booktabs}%
\usepackage{algorithm}%
\usepackage{algorithmicx}%
\usepackage{algpseudocode}%
\usepackage{listings}%
\usepackage{cleveref}
\usepackage{tikz}

\newtheorem{theorem}{Theorem}[section]
\newtheorem{definition}[theorem]{Definition}
\newtheorem{lemma}[theorem]{Lemma}
 
\newtheorem{corollary}[theorem]{Corollary}
\newtheorem{example}[theorem]{Example}
\newtheorem{proposition}[theorem]{Proposition}
\newtheorem{observation}[theorem]{Observation}

\begin{document}

\title[Adjacency Spectrum and Wiener Index of Essential Ideal Graph of a Finite Commutative Ring]{Adjacency Spectrum and Wiener Index of Essential Ideal Graph of a Finite Commutative Ring}


\author[1]{\fnm{Jamsheena} \sur{P}}\email{panikkarajamsheena@gmail.com}

\author*[1]{\fnm{Chithra} \sur{A V}}\email{chithra@nitc.ac.in}

\author[2]{\fnm{Subarsha } \sur{Banerjee}}\email{subarshabnrj@gmail.com}

\affil[1]{\orgdiv{Department of Mathematics}, \orgname{National Institute of Technology Calicut}, \orgaddress{
\city{Kozhikode}, \postcode{673601}, \state{Kerala}, \country{India}}}

\affil[2]{\orgdiv{Department of Mathematics}, \orgname{ JIS University}, \orgaddress{
\city{Agarpara}, \postcode{700109}, \state{West Bengal}, \country{India}}}



\abstract{Let $R$ be a commutative ring  with unity. The essential ideal graph $\mathcal{E}_{R}$ of $R$, is a graph with a vertex set consisting of all nonzero proper ideals of \textit{R} and two vertices $I$ and $K$ are adjacent if and only if $I+ K$ is an essential ideal. In this paper, we  study the adjacency spectrum of the  essential ideal graph of the finite commutative ring $\mathbb{Z}_{n}$, for $n=\{p^{m}, p^{m_{1}}q^{m_{2}}\}$, where $p,q$ are distinct primes, and $m,m_{1}, m_2\in \mathbb N$. We show that $0$ is an eigenvalue of the adjacency matrix of $\mathcal{E}_{\mathbb{Z}_{n}}$ if and only if either $n= p^2$ or $n$ is not a product of  distinct primes. We also determine all the eigenvalues of  the adjacency matrix of $\mathcal{E}_{\mathbb{Z}_{n}}$ whenever  $n$ is a product of three or four distinct primes. Moreover, we calculate the topological indices, namely the Wiener index and hyper-Wiener index of the essential ideal graph of $\mathbb{Z}_{n}$ for different forms of $n$.}



\maketitle

\section{Introduction}\label{sec1}

In recent decades, researchers have been exploring algebraic structures using graph theory properties. The notion of a graph connected to the zero divisors of a commutative ring was put forward by Beck \cite{bec} in $1998$.
However, the current definition and term for the \textit{zero-divisor graph} were initially presented by Anderson and Livingston \cite{liv} in 1999.
Following this, various studies were pursued on graphs defined on commutative rings by taking the ideals as vertices. 
Another graph, namely the \textit{comaximal ideal graph} was introduced in \cite{Ye} as a graph with vertices as the proper ideals of the ring $R$, and a pair of vertices $I$ and $K$ are adjacent if and only if $I+K =R$. Interested readers may refer to the papers \cite{akbari2017co, beh1, sharma} for more details. In $2018$, J. Amjadi \cite{amj} introduced the \textit{essential ideal graph} of a commutative ring, which is denoted as $\mathcal{E}_R$. 
The graph $\mathcal{E}_R$ has all nonzero proper ideals of $R$ as the vertex set,  and any two vertices $I$ and $K$ are adjacent if and only if $I+K$ is an essential ideal.  It is worth mentioning that a proper ideal $I$ of a ring $R$ is said to be an \textit{essential ideal} if it  has a nonzero intersection with every other non-zero ideal of $R$. 

 In mathematical chemistry, molecular descriptors like the topological indices play a vital role. A topological index is an invariant of molecular graphs that can be used to study the properties of their constituent molecules.
 Among this, the \textit{Wiener index}  introduced by H. Wiener \cite{wnr}is a well-known molecular descriptor,  which, in particular, is used for the preliminary testing of drug molecules.
 A generalization of the Wiener index known as the \textit{hyper-Wiener index} was introduced by M. Randić \cite{randc}, and is widely used in biochemistry. 
 Determining various topological indices of graphs associated with different algebraic structures has been an interesting area of research in the past few years. To get a better understanding of this, refer \cite{banerj1, selv}. Being motivated by the previous works, in this paper, we determine the Wiener and the hyper-Wiener index of $\mathcal{E}_{\mathbb Z_n}$, where $\mathbb Z_n$ is the ring of integers modulo $n$.

Let $G$ be a simple graph with vertex set  $V(G)=\{v_1, v_2,\cdots,v_n\}$ and edge set $E(G)$. 
The \textit{degree} of a vertex $v$, denoted by $deg(v)$, is defined as  the number of edges that occur in $v$. 
In $G$, a vertex $v$ is said to be universal, if it is  adjacent to all other vertices. The  \textit{complete graph} $K_{n}$, is a graph in which any two vertices are adjacent.
A graph $G$ is a $k-\textit{partite graph}$ if $V(G)$ can be partitioned into $k$ subsets $V_{1},\ V_{2},\ \cdots,\ V_{k}$ (named partite sets) such that the vertices $u$ and $v$ form an edge in $G$ if they belong to different partite sets.
If, in addition, there exists an edge between every two vertices belonging to different partite sets, then graph $G$ can be classified as $\textit{complete k-partite graph}$. The graph denoted as $K_{m,n}$ represents a complete bipartite graph consisting of two sets with sizes $m$ and $n$ respectively.
The \textit{induced subgraph}, $G[S]$, is formed by taking the subset $S$ of vertices from $G$, along with all the edges that connect vertices solely within $S$.
The \textit{complement} of a graph $G$ is denoted by  $\overline{G}$.
A set of vertices in a graph $G$ is \textit{independent} if any two vertices in the set are nonadjacent.
The \textit{join} of two graphs, $G_{1}$ and $G_{2}$, represented as $G_{1}\vee G_{2}$,  is formed by adding edges between any two vertices $v_1$ and $v_2$, where  $v_1\in G_{1}$ and  $v_2\in G_{2}$.
The \textit{adjacency matrix} $A(G)$ of a graph $G$ of order $n$  is the $n\times n$ matrix $A(G)= (a_{ij})$, where  $a_{i,j} = 1 $, if $v_i$ is adjacent to  $v_j$ in  $G $ and $a_{ij}= 0$ otherwise. The \textit{eigenvalues} of a graph $G$ are defined to be the eigenvalues of its adjacency matrix. The collection of all eigenvalues of $G$ is called the \textit{spectrum(adjacency spectrum) } of $G$. The \textit{energy} of a graph $G$, denoted by $\mathbb{E}(G)$, is defined as the sum of the absolute values of the eigenvalues of $A(G)$.
That is,  $\mathbb{E}(G)=\displaystyle \sum_{i=1}^{n}\lvert \lambda_i \rvert$, where $\lambda_1, \lambda_2, \cdots, \lambda_n$ are the eigenvalues of $G$. A graph  whose energy is greater than that of a complete graph is called \textit{hyperenergetic}. Otherwise, it is called \textit{non-hyperenergetic}.  
To delve into more definitions and results in the fields of ring theory, 
graph theory, and spectral graph theory,  one can refer  \cite{ ati, biggs, Ckovic, dummit, west, cha}. \\

  The paper is organized as follows. In Section $2$, we state the results that are needed for the subsequent sections. In Section $3$, we find the eigenvalues of $\mathcal{E}_{\mathbb{Z}_{n}}$ for $n=p^{m},m>1$, and  $n= p^{m_{1}}q^{m_{2}}$, where $p$ and $q$ are distinct primes with $p< q$ and $m_{1}, m_2$ are positive integers. We also prove that for the essential ideal graph of $\mathbb{Z}_{n}$, $0$ is not an eigenvalue if and only if either $n=p^m$, $m>2$ or $n$ is a product of distinct primes. In Section $4$, we calculate the Wiener index and the hyper-Wiener index of the essential ideal graph of $\mathbb{Z}_{n}$ for different values of $n$.
 
\section{Preliminaries}\label{sec2}
 The results shown below are beneficial for the parts that follow.
\begin{lemma}\cite{Zhang} \label{detMATX}
	Let $M,N,P,Q$ be matrices and let $Q$ be non-singular. Let $S= \begin{pmatrix}
	 M & N \\
	 P & Q  \\
	\end{pmatrix}$, then $\det S =  \det Q \times \det (M- NQ^{-1} P)$. Here, $(M- NQ^{-1} P)$ is known as  the \textbf{Schur complement} of $M$ in $S$.
\end{lemma}
\begin{lemma}
\cite{Zhang}
\label{block invrs}
Let $S \in M_{n}(\mathbb F)$ be a partitioned $2\times2$ block matrix $S= \begin{pmatrix}
	S_{11} & S_{12}  \\
	S_{21} & S_{22} \\
\end{pmatrix}$,where each matrix $S_{ii} \in M_{n_{i}}(\mathbb F)$, $i=1,2$ and $n= n_{1}+n_{2}$. If $S_{11},\ S_{22}$ and both Schur complements $S_{11}- S_{12}S_{22}^{-1} S_{21}$ and $S_{22}- S_{21}S_{11}^{-1}S_{12}$ are all invertible, then \\

$\begin{pmatrix}
	S_{11} & S_{12}  \\
	S_{21} & S_{22} \\
\end{pmatrix}^{-1} = \begin{pmatrix}
 (S_{11}- S_{12}S_{22}^{-1} S_{21})^{-1}  &  - S_{11}^{-1}S_{12}(S_{22}- S_{21}S_{11}^{-1}S_{12})^{-1} \\
 -(S_{22}- S_{21}S_{11}^{-1}S_{12})^{-1}S_{21}S_{11}^{-1} &  (S_{22}- S_{21}S_{11}^{-1}S_{12})^{-1} \\
\end{pmatrix}$.

\end{lemma}
\begin{proposition}\cite{magi} \label{circlnt det}
	Let $C_{(a,b,n)}= \begin{pmatrix}
    a & b & \cdots & \cdots & b\\
	b & a & b & \cdots  & b\\
	\vdots & \vdots & \vdots & \ddots & \vdots\\
	b & b & b & \cdots  & a \\
	\end{pmatrix}$ be a circulant matrix of order $n \times n$ with entries $a,b\in \mathbb{R}$. Then its determinant, denoted by $\delta$, is given by $\delta= (a+(n-1)b)(a-b)^{n-1}$.
 \end{proposition}
 \begin{proposition} \cite{magi} \label{crclnt invrs}
		If the circulant matrix $C_{(a,b,n)}$ is nonsingular, then its inverse is given by $C^{-1} _ {(a,b,n)}= \frac{1}{\delta} \begin{pmatrix}
		\delta_{n-1} & \Delta_{n-1} & \cdots& \Delta_{n-1} \\
		\Delta_{n-1}& 	\delta_{n-1} & \cdots& \Delta_{n-1} \\
		\vdots & \vdots & \ddots & \vdots \\
		\Delta_{n-1} & \cdots & \cdots & 	\delta_{n-1} \\
	\end{pmatrix}$, where 
 \\$\delta_{n-1}=(a+(n-2)b)(a-b)^{n-2}$ and $\Delta_{n-1} = -b(a-b)^{n-2}$.
	\end{proposition}
\begin{theorem}\cite{Ckvc} \label{chr.pol.rglrG}
	If \textit{G} is a regular graph of degree $\textit{r}$ with \textit{n} vertices, then the characteristic polynomial of $\overline{G}$ is $P_{\overline{G}}(\lambda)= (-1)^n \frac{\lambda-n+r+1 }{\lambda+r+1}  P_G(-\lambda-1)$.
\end{theorem}

\begin{theorem} \cite{Ckvc} \label{ch.pl.join}
	 Let $G_{1}$ and $G_{2}$ be two graphs of order $n_{1}$ and $n_{2}$ respectively.
    Then the characteristic polynomial of the join of $G_{1}$ and $G_{2}$  is given by  
    \begin{equation*}
        \begin{split}
            P_{G_{1}  \vee G_{2}}(\lambda)&= (-1)^{n_2} P_{G_1}(\lambda) P_{\overline{G_2}}(-\lambda-1) + 
        \\
        &(-1)^{n_1} P_{G_2}(\lambda) P_{\overline{G_1}}(-\lambda-1)- (-1)^{n_{1}+n_2}P_{\overline{G_1}}(-\lambda-1)P_{\overline{G_2}}(-\lambda-1).
        \end{split}
    \end{equation*}
\end{theorem}

\begin{observation} \cite{amj} \label{obsrvns}Let $R$ be a commutative ring with nonzero unity. Then every proper essential ideal of \textit{R} is a universal vertex in $\mathcal{E}_R$.
\end{observation} 
For any composite integer $n> 1$, let $n=p_{1}^{\alpha_{1}}p_{2}^{\alpha_{2}} \cdots p_{k}^{\alpha_{k}}$, $(k,\ \alpha_{i})\in \mathbb{N} $,  $(k,\alpha_{1})\neq (1,1)$, $p_{i}$'s are distinct primes ($1\le i\le k$). \\
\begin{theorem} \label{Strctr thm Zn} \cite{jam1}
For the essential ideal graph $\mathcal{E}_{\mathbb{Z}_{n}}$,
  $\mathcal{E}_{\mathbb{Z}_{n}}\cong H \vee K_{m}$, where $H$ is a $k$-partite graph and $K_{m}$ is a complete graph of order $m= \displaystyle \prod_{i=1}^{k}\alpha_{i}-1$.
\end{theorem}

\section{Adjacency Spectrum of  Essential Ideal Graph  of $\mathbb{Z}_{n}$}\label{sec3}
In this section, we study the adjacency spectrum of the essential ideal graph of $\mathbb{Z}_{n}$. We obtain the spectrum for $n= p^{m_{1}}$ and  $n= p^{m_{1}}q^{m_{2}}$, where $p$ and $q$ are distinct primes with $p<q$, and $m_1, m_2$ are positive integers. We also determine the adjacency  spectrum of $\mathcal{E}_{\mathbb{Z}_{n}}$, when $n$ is a product of three distinct primes and  a product of four distinct primes. 
Throughout the section, by spectrum of $\mathcal{E}(\mathbb Z_n)$, we shall mean the adjacency spectrum of $\mathcal{E}(\mathbb Z_n)$.

We  first prove the following result which provides a necessary and sufficient condition for an ideal to be an essential ideal of $\mathbb Z_n$.\\

\begin{theorem}  \label{chr. essentl.idl Zn}
   Let $n= p_{1}^{m_1}p_{2}^{m_2} \cdots p_{k}^{m_k}$ where $p_{1}<p_{2}< \cdots < p_{k}$ are distinct primes, and $m_i$ is a non-negative integer for $1 \le i \le k$. Any nonzero ideal $I= \langle p_{1}^{r_1}p_{2}^{r_2} \cdots p_{k}^{r_k} \rangle $ of $\mathbb {Z}_n$ is essential if and only if $r_i \ne m_{i}$ for any $i$. 
\end{theorem}
\begin{proof}
    Assume that  $I= \langle p_{1}^{r_1}p_{2}^{r_2} \cdots p_{k}^{r_k} \rangle$ be a nonzero essential ideal of $\mathbb {Z}_n$. We need to prove that $r_i \ne m_{i}$ for any $i$, $1 \le i \le k$. Suppose that  $r_i = m_{i} $ for some $i$; say $1$. Then for the ideal $N = \langle p_{2}^{m_2} p_{3}^{m_3} \cdots p_{k}^{m_k} \rangle$, $I\cap N = \langle 0 \rangle$, contradicting the fact that $I$ is essential. \\
  
  Conversely, let   $r_i \ne m_{i}$ for any $i$. That is, $I= \langle p_{1}^{r_1}p_{2}^{r_2} \cdots p_{k}^{r_k} \rangle$, $0 \le r_i \le m_{i}-1$ for $1 \le i \le k$ be a nonzero ideal. We need to prove that $I$ is essential. If not, there is a nonzero ideal $L \ne I$ such that $I\cap L= \langle 0 \rangle$. But  all the ideals of $\mathbb {Z}_n$ other than $I$ will be in any one of the following sets.
    \begin{equation*}
        \begin{split}
            A_{p_{1}^{m_1}}&=\{ \langle p_{1}^{m_1}p_{2}^{r_2} \cdots p_{k}^{r_k} \rangle ; 0 \le r_i \le m_i\ \text{ for } 2\le i \le k \}  
    \\
    A_{p_{2}^{m_2}}&=\{ \langle p_{1}^{r_1}p_{2}^{m_2} \cdots p_{k}^{r_k} \rangle ; 0 \le r_1 \le m_1 -1 \text{ and }  0 \le r_i \le m_i \    \text{ for } \  3 \le i \le k  \}.
        \end{split}
    \end{equation*}
     \\
In general, $A_{p_{j}^{m_j}}  =\{ \langle p_{1}^{r_1}p_{2}^{r_2} \cdots p_{j-1}^{r_{j-1}}p_{j}^{m_j}p_{j+1}^{r_{j+1}}\cdots p_{k}^{r_k} \rangle ; 
 0 \le r_i \le m_i -1, \ \text{ for } \ 1 \le i \le j-1  \ \text{ and} \   0 \le r_i \le m_i \ \text{ for } \  j+1 \le i \le k \}$; $1 \le j \le k$. Thus, $L$ must be in any of the sets  $A_{p_{1}^{m_1}},\ A_{p_{2}^{m_2}}, \cdots, A_{p_{k}^{m_k}} $ so that its intersection with $I$ is nonzero. This contradiction proves the result.
\end{proof}

\begin{proposition}\label{Aspctrm-p^m}
    	Let $n=p^m$, $m>2$ be a positive integer and $p$ be any prime. Then the spectrum of the essential ideal graph $\mathbb{Z}_{n}$ is $\begin{pmatrix}
	m-2 & -1 \\
	1 & m-2  \\
	\end{pmatrix}$.
\end{proposition}
\begin{proof}
If $n= p^m$, then all the nonzero proper ideals of $\mathbb{Z}_{n}$ are essential by Lemma \ref{chr. essentl.idl Zn},  and hence $\mathcal{E}_{\mathbb{Z}_{n}}$is a complete graph.
\end{proof}
\begin{corollary}
Let $n=p^m$, $m>2$ be a positive integer and $p$ be any prime. Then the energy of the essential ideal graph  of ${\mathbb{Z}_{n}}$ is $2m-4$. 
\end{corollary}
\begin{theorem}
    Let $n= p^{m_{1}}q^{m_{2}}$, where $p$ and $q$ are distinct primes with $p<q$ and $m_1, m_2$ are positive integers. Then the characteristic polynomial of $\mathcal{E}_{\mathbb{Z}_{n}} $ is given by  $P_{\mathcal{E}_{\mathbb{Z}_{n}}} (\lambda) = \lambda^{m_1 +m_2 -2}(\lambda+1)^{m_1 m_2 -2}P(\lambda)$, where $P(\lambda)= \lambda ^3+ (2-m_1m_2)\lambda^2 +[(1-m_1 m_2)(m_1+m_2)-m_1m_2]\lambda-{m_1}^2 {m_2}^2. $
\end{theorem}

\begin{proof}
 We can partition the vertex of $\mathcal{E}_{\mathbb{Z}_{n}}$ as follows: 
 
   $X= \{\langle p^rq^s \rangle: 0 \le r \le m_1 -1 ,\ 0\le s \le m_2 -1$ and $ (r,s)\ne (0,0)\}$\\ 
  $\hspace*{0.7cm}V_1 = \{ \langle p^{m_1}q^{s} \rangle: 0\le s \le m_2 -1 \}$ and\\
  $\hspace*{0.7cm} V_2 = \{ \langle p^{r}q^{m_2}\rangle : 0 \le r \le m_1 -1 \}$ so that  $V(\mathcal{E}_{\mathbb{Z}_{n}}) \simeq X\cup V_1 \cup V_2$.\\
  By Theorem \ref{chr. essentl.idl Zn} and Observation \ref{obsrvns}, $\mathcal{E}_{\mathbb{Z}_{n}}[X] \simeq K_{m_{1}m_{2}-1}$. 
   And, since $V_1$ and $V_2$ consist of independent vertices, $\mathcal{E}_{\mathbb{Z}_{n}}[V_1, V_2] \simeq K_{m_2,m_1}$
   Thus, by Theorem \ref{Strctr thm Zn}, $\mathcal{E}_{\mathbb{Z}_{n}} \simeq K_{m_1 m_2 -1} \vee K_{m_2, m_1} $.\\
    To find the characteristic polynomial of $\mathcal{E}_{\mathbb{Z}_{n}}$, we take $G_1= K_{m_1 m_2 -1}$ and $G_2 =  K_{m_2, m_1} $.  Then using Theorems \ref{chr.pol.rglrG} and \ref{ch.pl.join}, we have,
  \begin{equation*}
     \begin{split}
         P_{\mathcal{E}_{\mathbb{Z}_{n}}} (\lambda)&= (-1)^{m_1 +m_2} P_{G_1}(\lambda) P_{\overline{G_2}}(-\lambda-1) +
         (-1)^{m_1 m_2 -1} P_{G_2}(\lambda) P_{\overline{G_1}}(-\lambda-1)
         \\&- (-1)^{m_1 m_2 -1+ m_1 +m_2} P_{\overline{G_1}}(-\lambda-1)P_{\overline{G_2}}(-\lambda-1).
     \end{split}
 \end{equation*}  Here,   $\overline{G_1}$ is  the empty graph consisting of $m_1 m_2 -1$ vertices and $\overline{G_2}= K_{m_2} \cup K_{m_1}$. Hence,
 \begin{equation*}
     \begin{split}
         P_{\overline{G_1}}(-\lambda-1)&= (-1)^{m_1 m_2 -1} (\lambda+1)^ {m_1 m_2 -1}, \\
	P_{\overline{G_2}}(-\lambda-1)&= (-1)^{m_1 +m_2-2}(\lambda+m_1)(\lambda+m_2)\lambda^{m_1 +m_2-2}.
     \end{split}
 \end{equation*}
 Thus we have, 
 \begin{equation*}
     \begin{split}
         P_{\mathcal{E}_{\mathbb{Z}_{n}}} (\lambda)&= (-1)^{m_1 +m_2}(\lambda-m_1m_2+2)(\lambda+1)^{m_1m_2-2}(-1)^{m_1 +m_2-2}(\lambda+m_1) 
         \\
	&\times (\lambda+m_2)\lambda^{m_1 +m_2-2}+(-1)^{m_1 m_2 -1}\lambda^{m_1 +m_2-2}(\lambda^2 - m_1 m_2)(-1)^{m_1 m_2 -1}
        \\
	&\times (\lambda+1)^{m_1 m_2 -1}-(-1)^{m_1 m_2 -1+m_1 +m_2}(-1)^{m_1 m_2 -1}(\lambda+1)^{m_1 m_2 -1}
        \\
	&\times(-1)^{m_1 +m_2-2}(\lambda+m_1)(\lambda+m_2)\lambda^{m_1 +m_2-2}
        \\
	&= \lambda^{m_1 +m_2 -2}(\lambda+m_1)(\lambda+m_2)(\lambda+1)^{m_1 m_2 -2} (\lambda-m_1m_2+2)+\lambda^{m_1 +m_2 -2}
        \\
        &\times (\lambda^2 - m_1 m_2)(\lambda+1)^{m_1 m_2 -1}- \lambda^{m_1 +m_2 -2}(\lambda+1)^{m_1 m_2 -1}(\lambda+m_1)(\lambda+m_2).
     \end{split}
 \end{equation*}
	On simplifying, we obtain  
 \begin{equation*}
     P_{\mathcal{E}_{\mathbb{Z}_{n}}}(\lambda)= \lambda^{m_1 +m_2 -2}(\lambda+1)^{m_1 m_2 -2}P(\lambda),
 \end{equation*}
 where
 \begin{equation*}
     P(\lambda)= \lambda ^3+ (2-m_1m_2)\lambda^2 +[(1-m_1 m_2)(m_1+m_2)-m_1m_2]\lambda-{m_1}^2 {m_2}^2.
 \end{equation*}
\end{proof}
\begin{corollary}
Let $n= p^mq^m$,where $p$ and $q$ are distinct primes and $m>1$.  Then the spectrum of $\mathcal{E}(\mathbb{Z}_{n})$ is \\
	\begin{equation*}
 \begin{pmatrix}
	\frac{k+\sqrt{k^2+4m^3}}{2}	 & 0 & -1 & \frac{k-\sqrt{k^2+4m^3}}{2} & -m \\
	1 & 2m-2 & m^{2}-2 & 1 & 1 \\
	\end{pmatrix}, \text{where}\  k=(m^2+m-2).   
	\end{equation*}
\end{corollary}
\begin{example}
Let $n=36= 2^{2}3^{2}$. The vertex set  of  $\mathcal{E}_{\mathbb{Z}_{36}}$ (see Fig.$1$) is  \\
$V=\{\langle 2 \rangle, \langle 3 \rangle,\langle 4 \rangle,\langle 6 \rangle, \langle 9 \rangle, \langle 12 \rangle, \langle 18 \rangle\}$.
It can be partitioned as $V= X\cup X_1\cup X_2$, where
\begin{equation*}
    X= \{\langle 2 \rangle, \langle 3 \rangle,\langle 6 \rangle\}, V_1= \{\langle 4 \rangle,\langle 12 \rangle)\}, V_2= \{\langle 9 \rangle, \langle 18 \rangle\}. 
\end{equation*}
 Since $X$ contains all the proper essential ideals of $\mathbb{Z}_{36}$,  the subgraph $\mathcal{E}_{\mathbb{Z}_{36}}[X]$ is $K_3$. Now, $\mathcal{E}_{\mathbb{Z}_{36}}[V_1,\ V_2]=K_{2,2}$ and hence   $\mathcal{E}_{\mathbb{Z}_{36}}\simeq  K_3 \vee K_{2,2}$. 
Then the spectrum of $\mathcal{E}_{\mathbb{Z}_{36}}$ is
\begin{equation*}
   \begin{pmatrix}
	2+2\sqrt{3}	 & 0 & -1 & 2-2\sqrt{3} & -2 \\
	1 & 2 & 2 & 1 & 1 \\
	\end{pmatrix} 
\end{equation*}

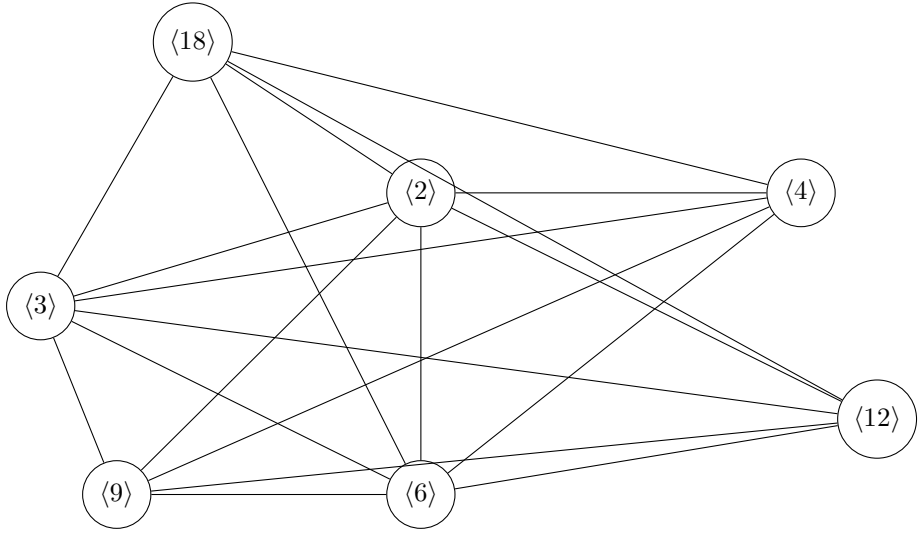
\begin{figure}[ht]
		\centering
		\begin{tikzpicture}
		\node[shape=circle,draw=black] (2) at (0,0) {$\langle 2 \rangle$};
		\node[shape=circle,draw=black] (3) at (-5,-1.5)  {$\langle 3 \rangle$};
		\node[shape=circle,draw=black] (6) at  (0,-4) {$\langle 6 \rangle$};
		\node[shape=circle,draw=black] (4) at  (5,0) {$\langle 4 \rangle$};
		\node[shape=circle,draw=black] (12) at (6,-3) {$\langle 12 \rangle$};
		\node[shape=circle,draw=black] (9) at (-4,-4) {$\langle 9 \rangle$};
		\node[shape=circle,draw=black] (18) at (-3,2) {$\langle 18 \rangle$};
		\draw (2) -- (3);
		\draw (6) -- (2);
		\draw (3) -- (6);

            \draw (2) -- (4);
            \draw (2) -- (12);
            \draw (2) -- (9);
            \draw (2) -- (18);

             \draw (3) -- (4);
            \draw (3) -- (12);
            \draw (3) -- (9);
            \draw (3) -- (18);

             \draw (6) -- (4);
            \draw (6) -- (12);
            \draw (6) -- (9);
            \draw (6) -- (18);

            \draw (12) -- (9);
            \draw (12) -- (18);
            \draw (18) -- (4);
            \draw (9) -- (4) ;

		\end{tikzpicture}
		\caption{$\mathcal{E}_{\mathbb{Z}_{36}}$}
		\label{Fig2}
		
	\end{figure}

\end{example}

\begin{corollary}
    Let $n= p^mq^m$,where $p$ and $q$ are distinct primes and $m>1$. Then the energy of the essential ideal graph of $\mathbb{Z}_{n}$ is $ k+\sqrt{k^2+4m^3}$, $k=(m^2+m-2)$.
\end{corollary}
\begin{lemma}\label{adjin prdct of primes}
Let $n= p_{1}p_{2} \cdots p_{k}$, where $p_{1}, p_{2}, \cdots, p_{k}$ are distinct primes. Then any two vertices $\langle x \rangle$ and $\langle y \rangle$ of the essential ideal graph of $\mathbb{Z}_{n}$ are adjacent if and only if $gcd(x,y)=1$, provided $x$ is the product of $i$ distinct primes  and $y$ is the product of $j$ distinct primes  for $1 \le i,j \le k-1 $.
\end{lemma}

\begin{theorem}
    Let $n= p_1p_2p_3$,  $p_i$ be a distinct prime for $1 \le i \le 3$. 
    Then the  spectrum of $\mathcal{E}_{\mathbb{Z}_{n}}$ is 
 \begin{equation*}
      \begin{pmatrix}
	1+\sqrt{2} &  \frac{-1+\sqrt{5}}{2} & 1-\sqrt{2} & \frac{-1-\sqrt{5}}{2}  \\
1 & 2 &1 & 2 \\
	\end{pmatrix}.
 \end{equation*}
\end{theorem}
	
\begin{proof}
    The number of nonzero proper ideals of $\mathbb{Z}_{n}$ is $2^3-2=6$. Hence, the adjacency matrix is a $6\times 6$ symmetric matrix.  Then the corresponding six vertices of $\mathcal{E}_{\mathbb{Z}_{n}}$ can be partitioned as follows. \\

$V_1=\{\langle p_1 \rangle, \langle p_2 \rangle,\langle p_3 \rangle\}$ and \\
$\hspace*{0.66cm}V_2= \{\langle p_2p_3 \rangle, \langle p_1p_3 \rangle,  \langle p_1p_2 \rangle\}$. 

By Lemma \ref{adjin prdct of primes}, we see that all vertices of $V_1$ are adjacent and they form the block  $J-I$ of order $3 $,  where $J$ is a matrix having all entries $1$, in the adjacency matrix of $\mathcal{E}_{\mathbb{Z}_{n}}$. Also, the vertices of $V_2$ are nonadjacent and each vertex is adjacent exactly to one of the vertices of $V_1$. Hence the vertices of $V_1$ and $V_2$ together form an identity  block  $I$  of order $3$ while the vertices of  $V_2$ form  a zero block of  order $3$. Then the adjacency matrix and characteristic polynomial of $\mathcal{E}_{\mathbb{Z}_{n}}$ are given by,

 \begin{equation*}
    A =\begin{pmatrix}
		J_{3\times 3}-I_{3\times 3} & I_{3\times 3} \\
			I_{3\times 3} &0_{3\times 3}\\
	\end{pmatrix}
 \end{equation*} 
 and  
 \begin{equation*} \det(A-\lambda I)
	= \det \begin{pmatrix}
	J-(\lambda +1)I & I \\
	I & -\lambda I \\
	\end{pmatrix}
  \end{equation*}
   \begin{equation*} 
	 =\det(J-(\lambda +1)I)\times \det(- \lambda I- I(J-(\lambda +1)I)^{-1}I),    \end{equation*} 
  where $J-(\lambda +1)I$ is the circulant matrix $C_{(- \lambda,1,3)}$.
By Propositions \ref{circlnt det} and \ref{crclnt invrs}, its  determinant and inverse are given by
 \begin{equation*}
     \delta=(2-\lambda)(\lambda+1)^{2}
 \end{equation*} and   
 \begin{equation*}
     \begin{split}
     C^{-1} _ {(-\lambda,1,3)}&=\frac{1}{(2-\lambda)(\lambda+1)^{2}}\begin{pmatrix}
	(\lambda ^2 - 1)& (\lambda +1) & (\lambda +1) \\
	 (\lambda +1) & (\lambda ^2 - 1)& (\lambda +1) \\
	  (\lambda +1) & (\lambda +1)& (\lambda ^2 - 1) \\
	\end{pmatrix} 
    \\
    &=\frac{1}{(2-\lambda)(\lambda+1)}\begin{pmatrix}
	(\lambda -1) & 1 & 1 \\
	1 & (\lambda -1) & 1 \\
	1 &  1 & (\lambda -1) \\
	\end{pmatrix}.
     \end{split}
 \end{equation*}
Also, 
 \begin{equation*}
     - \lambda I- I(J-(\lambda +1)I)^{-1}I= \begin{pmatrix}
     -\lambda -\frac{(\lambda-1)}{(2-\lambda)(\lambda+1)} & \frac{-1}{(2-\lambda)(\lambda +1)} & \frac{-1}{(2-\lambda)(\lambda +1)}\\
	\frac{-1}{(2-\lambda)(\lambda +1)} & 	-\lambda -\frac{(\lambda-1)}{(2-\lambda)(\lambda+1)} & \frac{-1}{(2-\lambda)(\lambda +1)}\\
\frac{-1}{(2-\lambda)(\lambda +1)} & \frac{-1}{(2-\lambda)(\lambda +1)} &	-\lambda -\frac{(\lambda-1)}{(2-\lambda)(\lambda+1)}\\
	\end{pmatrix},
 \end{equation*} is the circulant matrix $C_ {(	-\lambda -\frac{(\lambda-1)}{(2-\lambda)(\lambda+1)},  \frac{-1}{(2-\lambda)(\lambda +1)}, 3 )}$.\\
	Again by Proposition \ref{circlnt det},  
    \begin{equation*}
        \begin{split}
            \det(- \lambda I- (J-(\lambda +1)I)^{-1}) &= \frac{(-\lambda(2-\lambda)-1)(-\lambda(\lambda+1)+1)^{2}}{(2-\lambda)(\lambda+1)^{2}}, \text{ and  hence}
        \\
        \det (A-\lambda I)
	&=\frac{(2-\lambda)(\lambda+1)^{2}(\lambda^2-2\lambda-1)(\lambda^2+\lambda-1)^2}{(2-\lambda)(\lambda+1)^{2}}
 \\
 &= (\lambda^2-2\lambda-1)(\lambda^2+\lambda-1)^{2}.
        \end{split}
    \end{equation*}
  By solving the two quadratic polynomials, we obtain the required spectrum.
\end{proof} 
	 \begin{corollary}
	  Let $n= p_1p_2p_3$,  $p_i$ be a distinct prime for $1 \le i \le 3$. 
    Then
    \begin{enumerate}
        \item The energy of the essential ideal graph of $\mathbb{Z}_{n}$ is $2(\sqrt{2}+\sqrt{5})$.
        \item The graph $\mathcal{E}_{\mathbb{Z}_{n}}$ is non-hyperenergetic.	
    \end{enumerate}
    \end{corollary}

\begin{theorem}
Let $n= p_1p_2p_3p_4$, $p_i$ be a distinct prime for $1 \le i \le 4$.
Then the spectrum of the essential ideal graph of $\mathbb{Z}_{n}$ is \\
\begin{equation*}
    \begin{pmatrix}
	\frac{5+\sqrt{21}}{2} &  1 &   \frac{5-\sqrt{21}}{2} & \frac{-3+\sqrt{5}}{2} &  -1 &  \frac{-3-\sqrt{5}}{2}   \\
 1 & 5 & 1 & 3 & 1 & 3 \\
	\end{pmatrix}.
\end{equation*}
\end{theorem}	
\begin{proof}
To find the adjacency matrix, we first partition the vertex set of $\mathcal{E}_{\mathbb{Z}_{n}}$ as follows.\\		
	$V_	1 = \{ \langle p_1 \rangle, \langle p_2 \rangle,\langle p_3 \rangle, \langle p_4 \rangle \}	\\
	V_2 = \{\langle p_1p_2 \rangle, \langle p_1p_3 \rangle, \langle p_1p_4 \rangle, \langle p_2p_3 \rangle, \langle p_2p_4 \rangle, \langle  p_3p_4\rangle \}\\
	V_3 = \{ \langle p_1p_2p_3\rangle, \langle p_1p_2p_4\rangle, (p_1p_3p_4\rangle, \langle p_2p_3p_4\rangle \}$.\\	
	 Then by Lemma \ref{adjin prdct of primes}, we observe that all vertices of $V_1$ form the block matrix $J-I$ of order $4$ and each vertex of $V_2$ is adjacent exactly to one of the vertices of $V_2$, which will form the block matrix \\
 $  C = \begin{pmatrix}
	0 & 0 & 0 & 0 & 0 & 1 \\
	0 & 0 & 0 & 0 & 1 & 0\\
	0 & 0 & 0 & 1 & 0 & 0 \\
	0 & 0 & 1 & 0 & 0 & 0 \\
	0 & 1 & 0 & 0 & 0 & 0 \\
	1 & 0 & 0 & 0 & 0 & 0 \\
	\end{pmatrix}.$
  Also, each vertex of $V_1$ is adjacent exactly to three vertices of $V_2$ and one of the vertices of $V_3$, forming the blocks $B= \begin{pmatrix}
	0 & 0 & 0 & 1 & 1 & 1 \\
	0 & 1 & 1 & 0 & 0 & 1 \\
	1 & 0 & 1 & 0 & 1 & 0 \\
	1 & 1 & 0 & 1 & 0 & 0 \\
	\end{pmatrix}$ of order $4\times 6$ and $I$, identity matrix of order $4$, respectively. Finally, vertices of $V_2$  together with  vertices of $V_3$ form zero block of order $6 \times 4$ and vertices of $V_3$ constitute zero blocks of order $4$ to the adjacency matrix of $\mathcal{E}_{\mathbb{Z}_{n}}$. Hence the adjacency matrix is 
 
  \begin{equation*}
   A = \begin{pmatrix}
	(J-I)_{4\times4} & B_{4\times 6} & I_{4\times4} \\
	B_{6\times4}^{T} & C_{6\times6} & 0_{6\times4} \\
	I_{4\times4} & 0_{4\times 6} & 0_{4\times4}\\
	\end{pmatrix}.
  \end{equation*}
   
Then \\
\[A-\lambda I
	=  \left(
\begin{array}{cc|c}
     J-(\lambda+1)I & B & I \\
B^{T} & C-\lambda I & 0\\ \hline
 I & 0 & -\lambda I \\ 
\end{array}
\right) =\begin{pmatrix}
	M & N \\
	P & Q \\
	\end{pmatrix}.\]

 By Lemma \ref{detMATX}, if $\lambda\ne 0$, 
	\begin{equation*}\label{eqn1}
	\det (A-\lambda I)= \det Q \times \det (M- NQ^{-1} P).
	\end{equation*}\\
    Now,
	\begin{eqnarray} \label{eqn2}
	\det(M- NQ^{-1} P)=\det \begin{pmatrix}
	(J-(\lambda+1)I+\dfrac{1}{\lambda}I)_{4\times4} & B_{4\times6} \\
	B^{T} & D-\lambda I_{6\times6} \end{pmatrix}
	\end{eqnarray}
	\\ where $D-\lambda I_{6\times6}= \begin{pmatrix}
	-\lambda I_{3\times 3} & E_{3 \times3} \\
	E_{3 \times3} & -\lambda I_{3\times 3}\\
	\end{pmatrix}$, $E= \begin{pmatrix}
	0 & 0 & 1 \\
	0 & 1 & 0 \\
	1 & 0 & 0 \\
	\end{pmatrix}.$  
    \\  Applying Lemmas \ref{detMATX}, \ref{block invrs} and Proposition \ref{circlnt det}, we have 
    
    \begin{equation*}
        \det(D-\lambda I) =\det \begin{pmatrix}
	-\lambda I_{3\times 3} & E_{3 \times3} \\
	E_{3 \times3} & -\lambda I_{3\times 3}\\
	\end{pmatrix} 
	= (\lambda^{2}-1)^3
    \end{equation*}
  and   for $\lambda \ne \pm 1$, 
    \begin{equation*}
    \begin{split}
        (D- \lambda I) ^{-1}= \begin{pmatrix}
	\frac{\lambda}{1-\lambda^{2}}I & \frac{1}{1-\lambda^{2}} E \\
	\frac{1}{1-\lambda^{2}} E & \frac{\lambda}{1-\lambda^{2}}I \\
	\end{pmatrix} 
    \\
    \text { and }
    B(D- \lambda I)^{-1}B^{T}= \begin{pmatrix}
	\frac{3\lambda}{1-\lambda^{2}} & \frac{2+\lambda}{1-\lambda^{2}} & \frac{2+\lambda}{1-\lambda^{2}} & \frac{2+\lambda}{1-\lambda^{2}} \\
	\frac{2+\lambda}{1-\lambda^{2}} & \frac{3\lambda}{1-\lambda^{2}} & \frac{2+\lambda}{1-\lambda^{2}} & \frac{2+\lambda}{1-\lambda^{2}} \\
	\frac{2+\lambda}{1-\lambda^{2}} & \frac{2+\lambda}{1-\lambda^{2}} & \frac{3\lambda}{1-\lambda^{2}} & \frac{2+\lambda}{1-\lambda^{2}} \\
	\frac{2+\lambda}{1-\lambda^{2}} & \frac{2+\lambda}{1-\lambda^{2}} & \frac{2+\lambda}{1-\lambda^{2}} & \frac{3\lambda}{1-\lambda^{2}} \\
	\end{pmatrix}= C_ {(\frac{3\lambda}{1-\lambda^{2}}, \frac{2+\lambda}{1-\lambda^{2}},4)}. 
     \end{split}
    \end{equation*}	
And, Equation (\ref{eqn2}) is,	\begin{equation*}
	  \begin{split}
	        \det (M- NQ^{-1} P)&= (\lambda^{2}-1)^3 \times   \det C_{(\frac{\lambda^{4}-5\lambda^{2}+1}{\lambda(1-\lambda^{2})}, \frac{-\lambda^{2}-\lambda-1}{1-\lambda^{2}}, 4) }
	 \\
  &= (\lambda^{2}-1)^3 \ \bigg(\frac{\lambda^{4}-3\lambda^{3}-8\lambda^{2}-3\lambda+1}{\lambda (1-\lambda^{2})}\bigg) \ \bigg(\frac{\lambda^{4}+\lambda^{3}-4\lambda^{2}+\lambda+1}{\lambda (1-\lambda^{2})}\bigg)^3. 
	  \end{split}
    \end{equation*}
Hence, 
    \begin{equation*}
        \begin{split}
            \det (A-\lambda I)= \lambda^{4} \times  (\lambda^{2}-1)^3 \ \bigg(\frac{\lambda^{4}-3\lambda^{3}-8\lambda^{2}-3\lambda+1}{\lambda (1-\lambda^{2})}\bigg) \ \bigg(\frac{\lambda^{4}+\lambda^{3}-4\lambda^{2}+\lambda+1}{\lambda (1-\lambda^{2})}\bigg)^3 \\
	= (\lambda^{8}-9\lambda^{7}+26\lambda^{6}-29\lambda^{5}+29\lambda^{3}-26\lambda^{2}+9\lambda-1)(\lambda^{2}+3\lambda +1)^3.
        \end{split}
    \end{equation*} 
\end{proof}


\begin{corollary}
Let $n= p_1p_2p_3p_4$, $p_i$ be a distinct prime for $1 \le i \le 4$. Then 
\begin{enumerate}
    \item The energy of the essential ideal graph of $\mathbb{Z}_{n}$ is $20$.
    \item The graph $\mathcal{E}_{\mathbb{Z}_{n}}$ is non-hyperenergetic.	
\end{enumerate}
\end{corollary}

\begin{theorem}
\label{T1}
    Let $n= p_{1}^{m_1}p_{2}^{m_2} \cdots p_{k}^{m_k}$, where $p_{1}, p_{2}, \cdots, p_{k}$ are distinct primes and $m_i$ is a non negative integer for $1 \le i \le k$. Then $0$ is not an eigenvalue of $\mathcal{E}_{\mathbb{Z}_{n}}$ if and only if either $k=1$ and $m_1> 2$ or $m_i = 1$ for $1 \le i \le k$. 
\end{theorem}
\begin{proof}
    Assume that $0$ is not an eigenvalue of the adjacency matrix of $\mathcal{E}_{\mathbb{Z}_{n}}$.
    If  $n= p^m$, $m>2$, then we are done. Suppose that $n\ne p^m$, $m>2$. Then we need to prove that $m_i = 1$ for all $i=1,2, \cdots, k$. If possible, suppose $m_i > 1$ for atleast one $i$, say $m_1$. Without loss of generality, we assume that $n= p_{1}^{m_1}p_{2}\cdots p_{k}$. Then by Theorem \ref{chr. essentl.idl Zn}, the set of all essential ideals of $\mathbb{Z}_{n}$ is given by, \\ $X = \{ \langle p_1 \rangle, \langle p_{1}^{2}\rangle, \cdots, \langle p_{1}^{m_{1}-1}\rangle \}$. By Observation \ref{obsrvns}, these are the universal vertices of the essential ideal graph of $\mathbb{Z}_{n}$. Now consider the vertices $I= \langle p_{1}^{m_{1}-1}p_{2}\cdots p_{k}\rangle$ and $L= \langle p_{2} p_{3} \cdots p_{k}\rangle$. In $\mathcal{E}_{\mathbb{Z}_{n}}$, $I$ and $L$ are nonadjacent and are adjacent to any other vertex $K$ if and only if their sums $I+K$ and $L+K$ is an element of the set $X$.
    By elementary number theory, $K$ can be either an element of the set $X$ or the ideal $\langle p_{1}^{m_1}\rangle$.
   In other words, the adjacency and non-adjacency of the two vertices $I$ and $L$ are the same. Then the rows and columns  corresponding to the vertices $I$ and $L$ in the adjacency matrix are the same. Hence the matrix is singular and zero is an eigenvalue.
   
   Conversely, by Proposition \ref{Aspctrm-p^m}, the result is obvious when $n= p^m$, $m>2$.
    Now, let $n=p_1p_2p_3\cdots p_{k-1}p_k$.
   We shall index the rows and columns of the adjacency matrix of $\mathcal{E}_{\mathbb Z_n}$ in the following way:

    Let us consider the set $S =\{p_1,p_2,p_3,\ldots, p_{k-1},p_k\}$.
    Clearly, $S$ has $k$ elements.
   We first list the vertices of the form $\langle p_1\rangle, \langle p_2\rangle, \langle p_3\rangle, \ldots, \langle p_{k-1} \rangle, \langle p_k\rangle$. That is, we choose one element at a time from $S$.
   Next, we shall consider the vertices of the form $\langle p_1p_2\rangle,\langle p_1p_3\rangle,\cdots, \langle p_1p_k\rangle, \cdots, \langle p_{k-1}p_k\rangle$. That is, we shall choose two elements at a time from $S$.
   Clearly, we shall have $\binom{k}{2}$ such vertices.
   This process continues until we have exhausted all the vertices of $\mathcal{E}_{\mathbb Z_n}$.
Thus, in the end, we shall choose vertices of the form $\langle p_2p_3\cdots p_{k-1}p_k\rangle$, $\langle p_1p_3\cdots p_{k-1}p_{k}\rangle, \cdots$, $\langle p_1p_2\cdots p_{k-1}\rangle$. That is, we choose $k-1$ elements at a time from $S$ making a total of $\binom{k}{k-1}$ such vertices. Using the above indexing and Lemma \ref{adjin prdct of primes}, the adjacency matrix of $\mathcal{E}_{\mathbb Z_n}$ will be of the following form:
   \begin{equation}
         \label{E1}
         \left(
   \begin{array}{ccccccccccccccccccccccccccccc}
        & \dots & \dots & & \dots & & \dots & \dots & \dots & I_{k\times k}  \\
        & \dots & \dots & & \dots & & \dots & \dots & I_{\binom{k}{2}\times \binom{k}{2}} & 0\\
        & \dots & \dots & & \dots & & \dots & I_{\binom{k}{3}\times \binom{k}{3}} & 0 & 0\\
        & \dots & \dots & & \dots & & \dots & 0& 0& 0\\
        & \dots & \dots & & \dots & & 0 & \vdots& \vdots & \vdots\\
        & \dots & \dots & & \dots & & \vdots & \vdots & \vdots& \vdots\\
        & \dots & I_{\binom{k}{k-2}\times \binom{k}{k-2}} & &\vdots & & \vdots &
\vdots& \vdots & \vdots \\
        & I_{\binom{k}{k-1}\times \binom{k}{k-1}} & 0& & 0 & & 0& 0 & 0  & 0\\
   \end{array}
   \right).
   \end{equation}
Note that the matrix in  (\ref{E1}) is non-singular. Therefore, $0$ is not an eigenvalue of the adjacency matrix of $\mathcal{E}(\mathbb Z_n)$.
   This proves the result.
\end{proof}

\section{The Wiener and hyper-Wiener Index of the essential ideal graph of $\mathbb{Z}_{n}$ }
In this section, we compute the Wiener index and the hyper-Wiener index of $\mathcal{E}(\mathbb Z_n)$ for various $n$.
\begin{definition}
	The \textbf{Wiener index} of a graph $G$ is the sum of all distances between any pair of vertices of $G$. That is, \[ WW(G)= \sum_{u,v\in V(G)} d(u, v) = \frac{1}{2}\sum_ {u\in V(G)}d_G (u), \] where $d_G (u)$ is the sum of distances between $u$ and all other vertices of $V(G)$.
\end{definition}
\begin{definition}
The \textbf{ hyper-Wiener index} of a graph $G$ is defined as
\[ WW(G) = \frac{1}{2} W(G)+ \frac{1}{2} \sum_{u,v\in V(G)}  d^{2} (u,v).\]
\end{definition}


\begin{proposition}
    Let $n= p^{m}$, $m>1$ is a positive integer. Then 
    \begin{equation*}
    W(\mathcal{E}_{\mathbb{Z}_{n}})= WW(\mathcal{E}_{\mathbb{Z}_{n}})= \binom{m-1}{2}.    
    \end{equation*}
\end{proposition}
\begin{proof}
    By Lemma \ref{chr. essentl.idl Zn} the essential ideal graph $\mathcal{E}_{\mathbb{Z}_{n}}$ is complete if $n= p^{m}$ and hence 
    \[  W(\mathcal{E}_{\mathbb{Z}_{n}})= \frac{1}{2} \sum _{\langle x \rangle \in V(\mathcal{E}_{\mathbb{Z}_{n}})} d_{\mathcal{E}_{\mathbb{Z}_{n}}} (\langle x \rangle ), \]
    where \[d_{\mathcal{E}_{\mathbb{Z}_{n}}} (\langle x \rangle )= \sum_{\langle x \rangle \in V(\mathcal{E}_{\mathbb{Z}_{n}}) \ \langle y \rangle \ne \langle x \rangle } d(\langle x \rangle, \langle y \rangle)= m-2. \]
    Also,   \[ WW(\mathcal{E}_{\mathbb{Z}_{n}})= \frac{1}{2}\binom{m-1}{2} + \frac{1}{4} \sum_{ \langle x \rangle \in V(\mathcal{E}_{\mathbb{Z}_{n}})}d^{2}_{\mathcal{E}_{\mathbb{Z}_{n}}} (\langle x \rangle ), \] where  $d^{2}_{\mathcal{E}_{\mathbb{Z}_{n}}} (\langle x \rangle ) $ is  the sum of squares of distances between $\langle x \rangle$ and all other vertices of $\mathcal{E}_{\mathbb{Z}_{n}}$. Hence, $WW(\mathcal{E}_{\mathbb{Z}_{n}})=  \frac{1}{2}\binom{m-1}{2} + \frac{1}{4} (m-1)(m-2)= \binom{m-1}{2}$.
\end{proof}

\begin{theorem}\label{W-2case}
	Let $n= p^{m_{1}}q^{m_{2}}$, where $p<q$ are distinct primes and $m_1,\ m_2$ are positive integers.  Then the  Wiener index of the essential ideal graph of $\mathbb{Z}_{n}$ is
 \begin{equation*}
	    W(\mathcal{E}_{\mathbb{Z}_{n}}) = \frac{1}{2} [m_{1}m_{2} (m_1m_2 -1)+(m_{1}+m_{2})(2m_{1}m_{2}-4)+2(1+m_{1}^2 + m_{2}^2)].
	\end{equation*}
\end{theorem}
\begin{proof}
	First we partition the vertex set of $\mathcal{E}_{\mathbb{Z}_{n}} $ as follows :\\
 
	 $X= \{\langle p^rq^s \rangle: 0 \le r \le m_1 -1 ,\ 0\le s \le m_2 -1$ and $ (r,s)\ne (0,0)\}$\\ 
  $\hspace*{0.7cm}V_1 = \{ \langle p^{m_1}q^{s} \rangle: 0\le s \le m_2 -1 \}$ and\\
  $\hspace*{0.7cm} V_2 = \{ \langle p^{r}q^{m_2}\rangle : 0 \le r \le m_1 -1 \}$.\\
 Then $X$ is the set of all essential ideals of $\mathbb{Z}_{n}$ and induces a complete subgraph $K_{m_1 m_2 -1} $. The vertices of $V_1$ and $V_2$ are independent  and induce a complete bipartite graph $K_{m_2, m_1} $. Thus $\mathcal{E}_{\mathbb{Z}_{n}} \simeq K_{m_1 m_2 -1} \vee K_{m_2, m_1} $.
	For  every $\langle x \rangle \in X$,  the sum of the distances to any vertex  $\langle y \rangle \in V(\mathcal{E}_{\mathbb{Z}_{n}})$ can be obtained as
	$\sum_{\langle y \rangle \in X, \langle y \rangle \ne \langle x \rangle }d(\langle x\rangle,\langle y\rangle) + \sum_{\langle y \rangle \in V_1}d(\langle x\rangle,\langle y\rangle) + \sum_{\langle y \rangle \in V_2}d(\langle x\rangle,\langle y\rangle)$\\
	$\hspace*{1cm} = \sum_{\langle y \rangle \in X, y\ne x} 1+\sum _{\langle y \rangle \in V_1}1+ \sum_{\langle y \rangle \in V_2} 1$\\
	$\hspace*{1cm} = m_1 m_2 + m_1+ m_2-2 $.\\
	For  every $\langle x \rangle  \in V_1$,  the sum of the distances to any vertex  $\langle y \rangle \in V(\mathcal{E}_{\mathbb{Z}_{n}})$ can be determined as 	$\sum_{\langle y \rangle \in X}d(\langle x\rangle,\langle y\rangle)  + \sum_{\langle y \rangle \in V_1, \langle y \rangle \ne \langle x \rangle}d(\langle x\rangle,\langle y\rangle) + \sum_{\langle y \rangle \in V_2}d(\langle x\rangle,\langle y\rangle) $\\
	$\hspace*{1cm}= \sum_{\langle y \rangle \in X}1 +\sum_{\langle y \rangle \in V_1, y\ne x} 2+ \sum_{\langle y \rangle \in V_2}1$\\
	$\hspace*{1cm}=m_1 m_2+ m_1+ 2m_2-3$.\\	
	Similarly, for every $\langle x \rangle  \in V_2$, the sum of the distances to any vertex $\langle y \rangle \in V(\mathcal{E}_{\mathbb{Z}_{n}})$ is $m_1 m_2+2 m_1+ m_2-3$.	
	Hence the Wiener index of the graph $\mathcal{E}_{\mathbb{Z}_{n}}$ is given by $W(\mathcal{E}_{\mathbb{Z}_{n}})	= \frac{1}{2}[\sum_{\langle x \rangle  \in X} (m_1m_2+m_1+m_2-2)+ \sum_{\langle x \rangle \in V_1}( m_1 m_2+ m_1+ 2m_2-3) $\\ $\hspace*{1.5cm}+\sum_{\langle x \rangle  \in V_2} (m_1 m_2+2 m_1+ m_2-3)]$\\
	$\hspace*{1.5cm}= \frac{1}{2}[(m_1m_2-1)((m_1m_2+m_1+m_2-2))+ m_2(m_1 m_2+ m_1+ 2m_2-3) $\\ $  \hspace*{1.5cm} +m_1 (m_1 m_2+2 m_1+ m_2-3)]$\\ 
	$\hspace*{1.5cm}= \frac{1}{2} [m_{1}m_{2} (m_1m_2 -1)+(m_{1}+m_{2})(2m_{1}m_{2}-4)+2(1+m_{1}^2 + m_{2}^2)]$.
\end{proof}

\begin{corollary}
Let $n= p^mq^m$, where $p$ and $q$ are distinct primes and $m>1$.  Then the Wiener index of the essential ideal graph of $\mathbb{Z}_{n}$ is  \begin{equation*}
   W(\mathcal{E}_{\mathbb{Z}_{n}}) = \frac{m^4+4m^3+3m^2-8m+2}{2}.   
\end{equation*}	
\end{corollary}
\begin{theorem}
  Let $n= p^{m_{1}}q^{m_{2}}$, where $p<q$ are distinct primes and $m_1,\ m_2$ are positive integers.  Then the hyper-Wiener index of the essential ideal graph of $\mathbb{Z}_{n}$ is 	\begin{equation*}
      WW(\mathcal{E}_{\mathbb{Z}_{n}}) = \frac{1}{2} [m_{1}m_{2} (m_1m_2 -1)+(m_{1}+m_{2})(2m_{1}m_{2}-5)+3(m_{1}^2 + m_{2}^2)+2].
  \end{equation*}   
\end{theorem}

\begin{proof}
    By definition, 
    \begin{equation}\label{HW1}
      WW(\mathcal{E}_{\mathbb{Z}_{n}}) =   \frac{1}{2}W(\mathcal{E}_{\mathbb{Z}_{n}}) + \frac{1}{4} \sum_{ \langle x \rangle \in V(\mathcal{E}_{\mathbb{Z}_{n}})}d^{2}_{\mathcal{E}_{\mathbb{Z}_{n}}} (\langle x \rangle )     \end{equation} where \[ d^{2}_{\mathcal{E}_{\mathbb{Z}_{n}}} (\langle x \rangle )= \sum_{\langle y\rangle\in V(\mathcal{E}_{\mathbb{Z}_{n}})} d^{2} (\langle x\rangle, \langle y\rangle)\]    
    That is, the sum of squares of distances between the vertex $\langle x\rangle$ and all other vertices of $\mathcal{E}_{\mathbb{Z}_{n}}$. Now, take the same partition of $V(\mathcal{E}_{\mathbb{Z}_{n}})$ into $X\cup V_1 \cup V_2 $ described in the proof of Theorem \ref{W-2case}.\\
    Case $1$: $\langle x \rangle \in X$\\ 
Since  $ d(\langle x \rangle, \langle y \rangle)=1 $ for any  $ \langle y \rangle \in X,\ V_1 \ or \ V_2$ , $ d^{2}_{\mathcal{E}_{\mathbb{Z}_{n}}} (\langle x \rangle )= m_1 m_2 + m_1+ m_2-2 $
 Case $2$: $\langle x \rangle \in V_1$\\
 $ d^{2}(\langle x \rangle, \langle y \rangle)= \begin{cases}
     1, & \text{when} \  \langle y \rangle \in X  \ or \  V_2 \\
     4, & \text{when} \  \langle y \rangle \in V_1 \ \textit{and} \  \langle y \rangle \ne \langle x \rangle
 \end{cases}$.\\ 
 Thus  $ d^{2}_{\mathcal{E}_{\mathbb{Z}_{n}}} (\langle x \rangle )= m_1 m_2 + m_1+4m_2-5 $. \\
 Case $3$: $\langle x \rangle \in V_2$\\
Here,  $ d^{2}(\langle x \rangle, \langle y \rangle)= \begin{cases}
     1, & \text{when} \   \langle y \rangle \in X  \ or \  V_1 \\
     4, & \text{when} \   \langle y \rangle \in V_ 2\ \textit{and} \  \langle y \rangle \ne \langle x \rangle
 \end{cases}$. \\
 Then,  $ d^{2}_{\mathcal{E}_{\mathbb{Z}_{n}}} (\langle x \rangle )= m_1 m_2 + m_2+4m_1-5 $. \\
From Equation (\ref{HW1})  and  Theorem  \ref{W-2case}, we get the required result.
\end{proof}

Next, we calculate the Wiener and hyper-Wiener indices of $\mathcal{E}_{\mathbb{Z}_{n}}$ for $n= p_{1}p_{2} \cdots p_{k}$, where  $p_{1}< p_{2}< \cdots < p_{k}$  are distinct primes, using the idea of equitable partition of vertices.

\begin{definition}\cite{Hmrs}
 For a graph $G$, A partition of vertices   $V(G)= V_1\cup V_2 \cup \cdots \cup V_k$ is said to be an equitable partition if each vertex in $V_i$ has the same number of neighbors in $V_j$ for any $i,j \in \{1, 2, \cdots , k \}$.
\end{definition}

For this, consider the set $S =\{p_1,p_2,p_3,\ldots, p_{k-1},p_k\}$. Then  going through the same process as in the proof of Theorem \ref{T1},  we can exhaust all the vertices of $\mathcal{E}_{\mathbb{Z}_{n}}$. Also, we can see that the vertices of $\mathcal{E}_{\mathbb{Z}_{n}}$ can be partitioned into an equitable partition. That is, \\
$V_1=\{\langle p_i\rangle : 1\le i \le k \}$\\
$V_2=\{\langle p_ip_j\rangle : 1 \le i \le k-1 \textit{and} \  i+1 \le j \le k  \}$\\
$V_3=\{\langle p_ip_jp_l\rangle : 1 \le i \le k-2, \ i+1 \le j \le k-1 \textit{and} \  j+1 \le l \le k \}$\\
$\vdots $\\
$V_{(k-1)}=\{\langle p_1p_2p_3\cdots p_{k-1}\rangle , \langle p_1p_2p_3\cdots p_{k-2}p_k\rangle ,\cdots, \langle p_2p_3\cdots p_{k-1}p_k\rangle \}$.\\
Clearly $\vert V_1 \vert = \binom{k}{1}, \ |V_2|= \binom{k}{2}, \cdots, \ \text{and} \  |V_{(k-1)}|= \binom{k}{k-1}$.

By  Lemma \ref{adjin prdct of primes}, we can see that any vertex in $V_1$ has $\binom{k-1}{1}$ neighbors in $V_1$, $\binom{k-1}{2}$ neighbors in $V_2,\cdots, \binom{k-1}{k-1}$ neighbors in $V_{k-1}$. In general, any vertex of the set $V_t$ has $\binom{k-t}{1}$ neighbors in $V_1$, $\binom{k-t}{2}$ neighbors in $V_2,\cdots, \binom{k-t}{k-1}$ neighbors in $V_{k-1}$ respectively. Hence this  makes an equitable partition of $V(\mathcal{E}_{\mathbb{Z}_{n}})$ into sets $V_1$, $V_2$, $\cdots$, $V_{k-1}$.

\begin{theorem} \label{WI3}
Let $n= p_{1}p_{2} \cdots p_{k}$ where $p_{1}, p_{2}, \cdots, p_{k}$ are distinct primes. Then the Wiener index of the essential ideal graph $\mathcal{E}_{\mathbb{Z}_{n}}$ is  \[ W(\mathcal{E}_{\mathbb{Z}_{n}}) = \frac{1}{2} \sum_{t=1}^{k-1} \binom{k}{t} [2^{k+1} + 2^t- 2^{k-t}-7].\] 
\end{theorem}

\begin{proof}
By definition,   \[ W(\mathcal{E}_{\mathbb{Z}_{n}})= \frac{1}{2} \sum_{\langle x \rangle \in V(\mathcal{E}_{\mathbb{Z}_{n}})} d_ {\mathcal{E}_{\mathbb{Z}_{n}}} (\langle x \rangle), \] where $d_ {\mathcal{E}_{\mathbb{Z}_{n}}} (\langle x \rangle)$ is the sum of distances between the vertex  $\langle x \rangle$ and all other vertices of $\mathcal{E}_{\mathbb{Z}_{n}}$.

 Let $\langle x \rangle \in V_t$,  for  $1\le t \le k-1$. 
Then,
\begin{equation*}
\begin{split}
d_ {\mathcal{E}_{\mathbb{Z}_{n}}} (\langle x \rangle)=  \sum_{\langle y\rangle \in V_{1}}  d(\langle x\rangle,\langle y\rangle) + \sum_ {\langle y\rangle \in V_{2}}  d(\langle x\rangle,\langle y\rangle)+\cdots + \sum_ {\langle y\rangle \in V_{t}, \ \langle y\rangle \ne \langle x \rangle} d(\langle x\rangle,\langle y\rangle)
\\
\sum_ {\langle y\rangle \in V_{t+1}} d(\langle x\rangle,\langle y\rangle) + \cdots +\sum_ {(y)  \in V_{k-t}}  d(\langle x\rangle,\langle y\rangle)+ \sum_ {\langle y\rangle  \in V_{k-t+1}}  d(\langle x\rangle,\langle y\rangle)+
\end{split}
\end{equation*}
 \begin{eqnarray}\label{WEqn1}
 \cdots+ \sum_{\langle y\rangle  \in V_{k-1}}d(\langle x\rangle,\langle y\rangle) 
 \end{eqnarray}
 By Lemma \ref{adjin prdct of primes},
if $\langle y\rangle \in V_s$ for $1\le s\le k-t$  then, \\
$ d(\langle x\rangle,\langle y\rangle)= \begin{cases}
1 & if gcd(x,y)= 1 \\
2 & if gcd(x,y)\ne 1 
\end{cases} $  and if $\langle y\rangle  \in V_{k-t+s}$ for $1\le s\le t-1$ then, \\

$ d(\langle x\rangle,\langle y\rangle)= \begin{cases}
3 & if gcd(x,y)=  \textit{product of s distinct primes} \\
2 & otherwise 
\end{cases} $ .\\
Then, for $1\le s \le k-t$; $s\ne t$,
\[  \sum_ {\langle y\rangle \in V_{s}}d(\langle x\rangle,\langle y\rangle)= 2  \binom{k}{s} -  \binom{k-t}{s} \]  and  for $1\le s \le t-1$, \[  \sum_ {\langle y\rangle \in V_{k-t+s}}d(\langle x\rangle,\langle y\rangle) = 2  \binom{k}{k-t+s} + \binom{t}{s}. \] 
 Hence by Equation (\ref{WEqn1}),
\[d_ {\mathcal{E}_{\mathbb{Z}_{n}}} (\langle x\rangle)= \sum_ {s= 1 \\
{s \ne t}}^{k-t} [2 \binom{k}{s}- \binom{k-t}{s}] +[2 \binom{k}{t}- \binom{k-t}{t}-2] + \sum_{s=1}^{t-1}[2\binom{k}{k-t+s} +\binom{t}{s}] \]
\[= 2\sum_{s=1}^{k-1} \binom{k}{s} - \sum _{s=1}^{k-t} \binom{k-t}{s} +\sum_{s=1}^{t-1} \binom{t}{s}-2= 2^{k+1} + 2^t- 2^{k-t}-7 .\]
Hence,  \[ W(\mathcal{E}_{\mathbb{Z}_{n}}) = \frac{1}{2} [\sum_{\langle x\rangle \in V_{1}} d_ {\mathcal{E}_{\mathbb{Z}_{n}}} (\langle x\rangle) +  \sum_{\langle x\rangle\in V_{2}} d_ {\mathcal{E}_{\mathbb{Z}_{n}}} (\langle x\rangle) +\cdots+ \sum_{\langle x\rangle \in V_{k-1}} d_ {\mathcal{E}_{\mathbb{Z}_{n}}} (\langle x\rangle)] \]
\[  = \frac{1}{2} \sum_{t=1}^{k-1} \binom{k}{t} [2^{k+1} + 2^t- 2^{k-t}-7]. \]

\end{proof}
\begin{theorem}
   Let $n= p_{1}p_{2} \cdots p_{k}$ where $p_{1}< p_{2}< \cdots < p_{k}$ are distinct primes. Then the hyper-Wiener index of the essential ideal graph $\mathcal{E}_{\mathbb{Z}_{n}}$ is  \[ WW(\mathcal{E}_{\mathbb{Z}_{n}}) = \frac{1}{2} \sum_{t=1}^{k-1} \binom{k}{t} [3\times 2^{k} - 2\times2^{k-t} + 3 \times2^t-13]. \]    
\end{theorem}
\begin{proof}
    By definition, 
    \begin{equation}\label{HW2}
         WW(\mathcal{E}_{\mathbb{Z}_{n}}) = \frac{1}{2} W(\mathcal{E}_{\mathbb{Z}_{n}})+ \frac{1}{4} \sum_{\langle x\rangle \in V(\mathcal{E}_{\mathbb{Z}_{n}})}  d_ {\mathcal{E}_{\mathbb{Z}_{n}}}^{2}(\langle x\rangle) 
    \end{equation}
    where $ d_ {\mathcal{E}_{\mathbb{Z}_{n}}}^{2}(\langle x\rangle) $ is the sum of squares distances between $\langle x\rangle$ and all other vertices of $\mathcal{E}_{\mathbb{Z}_{n}}$. \\
Let  $\langle x\rangle$ be a vertex of  $V_t $, for $1 \le t \le k-1$. For  $\langle y\rangle \in V_s$, $1\le s\le k-t$, \\
\begin{equation*}
d^{2}(\langle x\rangle, \langle y\rangle)= \begin{cases}
1, & \textit{if} \ gcd(x,y)= 1 \\
4, & \textit{if} \ gcd(x,y)\ne 1 
\end{cases}     
\end{equation*} 
and for $\langle y\rangle \in V_{k-t+s}$, $1\le s\le t-1$,\\
\begin{equation*}
d^{2}(\langle x\rangle, \langle y\rangle) = \begin{cases}
9 & if \ gcd(x,y)= \textit{product of s distinct primes} \\
4 & \textit{otherwise} 
\end{cases} .    
\end{equation*} \\
Then,\\
\begin{eqnarray}\label{HW3}
    \begin{split}
     d_{\mathcal{E}_{\mathbb{Z}_{n}}} ^{2}(\langle x\rangle)= \sum_{\langle y\rangle\in V_1} d^{2} (\langle x\rangle, \langle y\rangle)  + \sum_ {\langle y\rangle \in V_{2}} d^{2} ((x), (y))+\cdots + \sum_ {\langle y\rangle \in V_{t}, \ (y)\ne (x)}  d^{2} (\langle x\rangle, \langle y\rangle) +
     \\
     \cdots + \sum_ {\langle y\rangle \in V_{k-t}}  d^{2}(\langle x\rangle, \langle y\rangle)  + \cdots+ \sum_ {\langle y\rangle\in V_{k-1}}  d^{2}(\langle x\rangle, \langle y\rangle)
    \end{split}
\end{eqnarray}

For $1\le s \le k-t$; $s\ne t$,  
\[  \sum_ {\langle y\rangle\in V_{s}} d^{2}(\langle x\rangle, \langle y\rangle) = 4  \binom{k}{s} -  3 \binom{k-t}{s} \] and  for $1\le s \le t-1$, \[  \sum_ {\langle y\rangle \in V_{k-t+s}} d^{2} (\langle x\rangle, \langle y\rangle)  = 4 \binom{k}{k-t+s} +5  \binom{t}{s}. \] 
Finally, for $\langle y\rangle \in V_t $, 
\[  \sum_ {\langle y\rangle \in V_{t}} d^{2}(\langle x\rangle, \langle y\rangle) = 4  \binom{k}{t} -  3 \binom{k-t}{t}-4. \] 
Hence by Equations  (\ref{HW2}) and (\ref{HW3}) and Theorem $\ref{WI3}$, we get the required result. 
\end{proof}

\section{Conclusion}
In this paper, the adjacency spectrum of the  essential ideal graph of finite commutative ring $\mathbb{Z}_{n}$, for $n=\{p^{m}, p^{m_{1}}q^{m_{2}}\}$, where $p,q$ are distinct primes, and $m,m_{1}, m_2\in \mathbb N$ is determined.  All the eigenvalues of $\mathcal{E}_{\mathbb{Z}_{n}}$ whenever $n$ is a product of three or four distinct primes are computed.
Further, we have established a characterization of the ring $\mathbb{Z}_{n}$ for which $0$ is an eigenvalue of $\mathcal{E}_{\mathbb{Z}_{n}}$.
In addition, the topological indices, namely the Wiener index and hyper-Wiener index of the essential ideal graph of $\mathbb{Z}_{n}$ for different forms of $n$ are calculated.

\section{Declarations}
\textbf{Conflict of interest} On behalf of all authors, the corresponding author states that there is no conflict of interest.

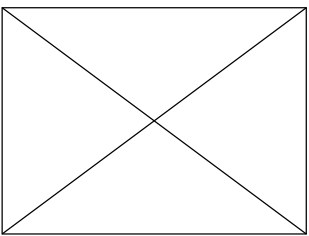
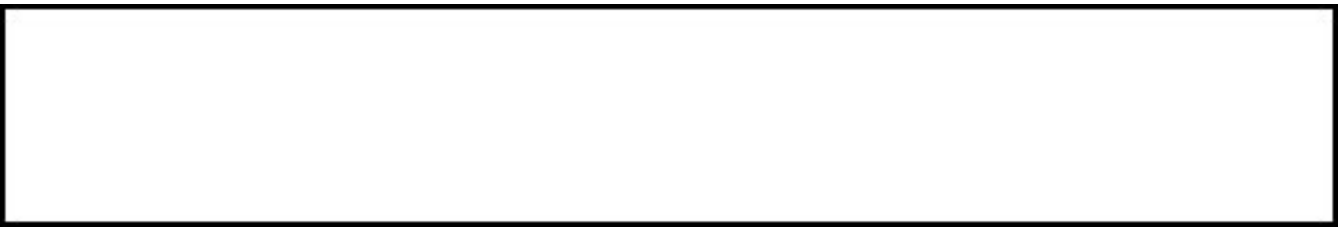
\end{document}